\documentclass[12pt,oneside]{amsart}
\usepackage{amsfonts}\usepackage{amssymb}\usepackage{amsmath}
\usepackage{hyperref}
\usepackage{epsfig}
\usepackage[dvips]{color}

\def\dim{{\rm dim\,}}
\def\Hom{{\mathop{\rm Hom}}}

\newcommand{\ch}{{\mathfrak{h}}}

\newcommand{\cn}{{\mathfrak{n}}}

\newcommand{\cz}{{\mathfrak{z}}}

\newcommand\eS{{\sf S}_+^2}
\newcommand\eeS{{\sf S}^2}
\def\g{\mathfrak{g}}
\def\so{\mathop{\mathfrak{so}}}\def\su{\mathop{\mathfrak{su}}\nolimits}
\def\N{\cn}
\def\LLL{\widehat{\mathcal{L}}}
\def\NNN{\widehat{\mathcal{N}}}
\def\NN{\mathcal{N}}
\def\DD{\mathcal{D}}
\def\GL{{\mathop{\rm GL}}}
\def\Z{\mathfrak{z}}
\def\H{\mathfrak{h}}
\def\O{{\mathop{\rm O}}}\def\SO{{\mathop{\rm SO}}}

\def\id{{\rm id}}

\def\det{{\rm det\,}}

\def\Gr{{\rm Gr}}

\def\Im{{\mathop{\rm Im}}}

\def\Spec{{\mathop{\rm Spec}}}
\def\Tr{{\mathop{\rm Tr}}}

\def\Tr{{\mathop{\rm Tr}}}

\def\R{\mathbb{R}}
\def\C{\mathbb{C}}

\def\RP{\mathbb{R}{\rm P}}

\def\eqref#1{(\ref{#1})}

\newtheorem{theorem}{Theorem}[section]

  \newtheorem{remark}{Remark}

\author{Sergio Console * , Anna Fino *, Evangelia Samiou **}
\thanks{* Research partially supported by MURST, GNSAGA (Indam) of 
Italy and EDGE Research Training Network
HPRN-CT-2000-00101}
\thanks{** Research  supported by  GNSAGA (INDAM) of Italy and the 
University of Cyprus}


\address{
Sergio Console, Anna Fino\\
Dipartimento di Matematica\\
Universit\`a di Torino \\
via Carlo Alberto, 10\\
10123 Torino, Italy}

\email{sergio.console@unito.it,  annamaria.fino@unito.it}

\address{Evangelia Samiou \\
University of Cyprus \\
Department of Mathematics and Statistics \\
P.O. Box 20537\\
1678 Nicosia, Cyprus}
\email {samiou@ucy.ac.cy}

\begin{document}

\title[The moduli space]{The moduli space of 6-dimensional 2-step nilpotent Lie
algebras}

\maketitle
\begin{abstract}
We determine the moduli space of metric $2$-step nilpotent Lie 
algebras of dimension up to $6$.
This space is homeomorphic to a cone over a $4$-dimensional 
contractible simplicial complex.
\end{abstract}


\section{Introduction}\label{Introduction}

The geometry of 2-step nilpotent Lie groups with a left-invariant metric is
very rich and has been widely studied since the papers of A. Kaplan \cite{k,
k2} and P. Eberlein \cite{e} (see, for instance \cite{df,df2, e2, e3, gms,
l2, sa2} for recent papers on this subject). Important examples of such
Lie groups are provided by groups of Heisenberg type \cite{k2}.

In general, the moduli space of metric Lie algebras of a fixed dimension
is a cone with peak the abelian Lie algebra and basis the subset
obtained by normalizing the Lie bracket $c$, for instance requiring 
$\Tr(c^* c)=2$.
In the present paper we determine the moduli space $\NN_6$ of
$6$-dimensional $2$-step nilpotent Lie algebras endowed with a metric.
We show that $\NN_6$ is a cone over an explicitely given contractible 
$4$-dimensional
simplicial complex. We also exhibit standard metric representatives of
the $7$ isomorphism types of $6$-dimensional $2$-step nilpotent Lie algebras
within our picture. This contains all deformations of these Lie
algebras, cf. \cite{l2}.

In \cite{l2} J. Lauret identified in a natural way each point of the 
variety of real Lie
algebras with a left-invariant Riemannian metric on a Lie group and
studied the interplay between invariant-theoretic and Riemannian
aspects of this variety. We show that on a certain subset of
$\NN_6$ the nullity of the Riemannian curvature tensor singles
out products.

The subspace $\NN_{n,k}\subset\NN_n$ of Lie algebras with 
$k$-dimensional commutator
ideal contains the subspace $\DD_{n,k}$ of algebras with isometric $c^*$
as a strong deformation retract.
For the algebras $\cn_{(\alpha_+, \alpha_-)}\in\DD_{6,2}$
we give the structure equations, write down the curvature tensor
and compute their infinitesimal rank, i.e. the minimal
nullity of the Jacobi operators. In \cite{sa2} it was proved that
groups of Heisenberg type have infinitesimal
rank one. We show that this is  also the case for any
$\cn_{(\alpha_+, \alpha_-)} \in \DD_{6,2}$ with the exception of
$\cn_{(1,1)}\cong {\mathfrak h}_3 \oplus {\mathfrak h}_3$ and
$\cn_{(1/2,1/2)} \cong \N_5\oplus\R$, both endowed with the product metric,
whose rank is two.

\section{Preliminaries}\label{Preliminaries}

A Lie algebra $\g$ is nilpotent if its central series ends, i.e.
in the sequence of ideals of $\g$ recursively defined by $\g^0:=\g$,
$\g^{i+1}:=[\g , \g^i]$ there is an integer $k$ such that $\g^k=0$. Then $\g$
is a $k$-step
nilpotent if $\g^k=0$ and $\g^{k-1}\neq 0$. Thus a $2$-step nilpotent
Lie algebra $\N$ is a Lie algebra such that its
commutator ideal $\N^1:=[\N , \N]$ is contained in its centre.

A left-invariant metric on a (simply connected) $2$-step nilpotent
Lie group $N$
is given by a scalar product $\langle \cdot, \cdot \rangle$ on its
Lie algebra $\N$. We will call such a Lie algebra ``metric $2$-step
nilpotent Lie algebra".

A simply connected $2$-step nilpotent Lie group
with left-invariant metric is uniquely determined by the triple
$(\H, \Z, j)$ \cite{e, k}, where $\H$ and $\Z$ are real
vector spaces with positive definite scalar product and $j\colon\Z\to\so (\H)$
is the homomorphism (of vector spaces, not necessarily of Lie
algebras) related to the
Lie bracket by
$$\langle [x,y], z \rangle=\langle y, j(z)x \rangle \ \ \ \ \forall
x,y\in\H, z\in\Z\ .$$
Thus, $j$ is essentially the adjoint of the Lie bracket
$c\colon\Lambda^2\H\to\Z $.
Requiring in addition $j$ to be injective makes this correspondence
one to one.

Observe that if one identifies $\cz$ with its dual $\cz^*$ via the
metric (and similarly for $\ch$ and $\N$), then the differential
$d: \Lambda^1\cz^*\subset \Lambda^1 \N^* \to \Lambda^2\ch^*\subset 
\Lambda^2 \N^*$
can be identified with $j$ and
$\dim (\Im \, d) = \dim ({\N}^1)$.

By \cite{M, GK} there are 34
classes of 6-dimensional nilpotent Lie algebras. Out of these 34
classes, the 2-step nilpotent have  the following structure
equations
$$
\begin{array} {l}
(0,0,0,12,13,23)\\
(0,0,0,0, 13 + 42, 14 + 23) = \ch_3^\C,\\
(0,0,0,0, 12, 14 + 23),\\
(0,0,0,0, 12, 34)=\ch_3\oplus\ch_3,\\
(0, 0,0,0, 12, 13)=\cn_5\oplus\R,\\
(0,0,0,0,0, 12 + 34)=\ch_5 \oplus \R,\\
(0,0,0,0,0, 12)=\ch_3\oplus\R^3,
\end{array}
$$
where  ${\mathfrak h}^{\C}_3$ is the complex 3-dimensional Heisenberg
Lie algebra, ${\mathfrak h}_3$  the real 3-dimensional
Heisenberg Lie algebra and
$\N_5 = (0,0,0,12,13)$. We use the notation of \cite{Sa}. For 
example, $(0,0,0,12,13,23)$ denotes
the Lie algebra with $de^i=0, i=1,\dots 3$, $de^4=e^1\wedge e^2$, 
$de^5=e^1\wedge e^3$,
$de^6=e^2\wedge e^3$, where $(e^j)$ is a basis of left-invariant 1-forms.

\section{The Moduli space of $2$-step nilpotent Lie
algebras}\label{The Moduli space of $2$-step
nilpotent Lie algebras}

We now determine the moduli space of metric $2$-step nilpotent Lie algebras.
Let $\g$ be a metric Lie algebra of dimension $n$. We can always
choose a linear isometry
of $\g$ with Euclidean space $\R^n$ (endowed with its standard scalar
product). The set
of Lie brackets on $\R^n$ is an algebraic subset of
$\Hom(\Lambda^2\R^n,\R^n)$, whose ideal is given by the Jacobi identity,
i.e.,
\[
\begin{array}{ll}
\LLL_n:=&\{c\in\Hom(\Lambda^2\R^n,\R^n)\mid           \\
            &c(c(u,v),w)+ c(c(w,u),v)+ c(c(v,w),u)=0\ \forall
u,v,w\in \R^n\} \ .
\end{array}
\]
The set of $2$-step nilpotent Lie brackets on $\R^n$ is
$$\NNN_n:=\{c\in\Hom(\Lambda^2\R^n,\R^n)\mid c(c(u,v),w)=0\ \forall\
u,v,w\in \R^n\} \ . $$
These sets are invariant under the $\GL(n, \R)$-action on
$\Hom(\Lambda^2\R^n,\R^n)$.
The moduli space of $2$-step nilpotent (metric)
$n$-dimensional Lie algebras is the space of (isometric) isomorphism
classes of such Lie algebras. It
inherits its topology as the quotient of $\NNN_n$ by the action of
$\GL(n, \R)$ (respectively, $\O(n)$),
$$\NN_n:=\NNN/\O(n) \qquad (\mbox{resp. }\check{\NN}_n = \NNN/\GL(n, 
\R) )  \  .   $$

For $k\leq n$ we decompose $\R^n=\R^{n-k}\oplus\R^{k}$ orthogonally.
A metric $n$-dimensional $2$-step nilpotent Lie algebra $\N$ with
$\dim \N^1 = k$ is isometric to $(\R^n,c)$ where $c$ is a
$2$-step nilpotent Lie bracket of rank $k=\dim(\Im c)$ such that $\Im c
= \{ 0 \}\oplus\R^k $.
We define
$$ \NNN_{n,k}:=\{c\in\Hom(\Lambda^2\R^{n-k},\R^k)\mid c \mbox{
surjective} \}.$$
This space carries an action of $\GL(n-k, \R)\times\GL(k, \R)$. The moduli
space of $n$-dimensional
(metric) $2$-step nilpotent Lie algebras with $k$-dimensional
commutator ideal is the
quotient
$$ \begin{array} {l}\NN_{n,k}=\NNN_{n,k}/(\O(n-k)\times\O(k))\\
(\mbox{resp. }
\check{\NN}_{n,k} = \NNN_{n,k}/(\GL(n-k, \R)\times\GL(k, \R)) ).
\end{array}$$
Extending $c\in\NNN_{n,k}$ by $0$ to all of $\Lambda^2\R^n$, we may view
$\NNN_{n,k}\subset\NNN_n$
and  decompose
\begin{equation}\label{decomp}
\NN_n = \bigcup_{ 0\leq k \leq \binom{n-k}{2} } \NN_{n,k} \qquad (\mbox{resp. }
\check{\NN}_n = \bigcup_{ 0\leq k \leq \binom{n-k}{2} } \check{\NN}_{n,k})  \ .
\end{equation}

We denote by $\gamma_{k,V} \to \Gr_k(V)$ the tautological vector bundle over
the Grassmanian of $k$-planes in a real vector space $V$. We let
$\eS\gamma_{k,V}^*\subset \eeS\gamma_{k,V}^* $ be the set of
positive definite symmetric
$2$-tensors on $ \gamma_{k,V}$.

The adjoint $c^*$ of $c\in\NNN_{n,k}$ is injective on $\R^k$. Pushing
forward the standard scalar product $g_\mathrm{std}$ on
$\R^k$,
    a scalar product on its image is defined.
The maps
\begin{equation}\label{maps}
\begin{array}{ccccc}
\NNN_{n,k}&\stackrel{\phi}{\longrightarrow}&\eS\gamma_{k,\Lambda^2\R^{n-k}}^*
&\stackrel{\pi}{\longrightarrow}&
\Gr_k(\Lambda^2\R^{n-k})  \\
    c                              &  \mapsto    &          (\Im c^*,
c^*g_\mathrm{std})       & \mapsto &  \Im c^*
\end{array}
\end{equation}
are $\O(n-k)\times\O(k)$-equivariant. In particular
\begin{theorem}\label{general}
There is a homeomorphism $$\NN_{n,k} \approx
\eS\gamma_{k,\Lambda^2\R^{n-k}}^*/\O(n-k)$$ and a
strong deformation retraction $$\NN_{n,k} \simeq
\Gr_k(\Lambda^2\R^{n-k})/\O(n-k)=:\DD_{n,k}.$$ Here 
$\DD_{n,k}\hookrightarrow\NN_{n,k}$
is identified with the subset of those $2$-step nilpotent Lie 
algebras with isometric
$j=c^*\colon\R^k\to\so(n-k)$.
\end{theorem}
\begin{proof} Two Lie brackets
$c,c^\prime\in\NNN_{n,k}\subset\Hom(\Lambda^2\R^{n-k},\R^k)$ are
equivalent in
$\NN_{n,k}$ if there are $A\in\O(k)$ and $T\in\O(n-k)$ such that
$AcT^{-1}=c^\prime $. Equivalent formulations are
$$
\begin{array} {c}
(T^{-1})^*c^*A^*= {c^\prime}^*,\\
c^*A^*\left( {c^\prime}^*|_{\Im {c^\prime}^*} \right)^{-1}=T^*,
\end{array}
$$
that is to say,
$T^*\colon\Im{c^\prime}^*\to\Im c^*$ is isometric with respect to the
metrics pushed forward by $c,c^\prime$.
Thus, the map
$\NNN_{n,k}\stackrel{\phi}{\longrightarrow}\eS\gamma_{k,\Lambda^2\R^{n-k}}^*$
in \eqref{maps} induces an homeomorphism on both quotients by
$\O(n-k)\times\O(k)$.

Let $g_0 $ be an
$\O(n-k)$-invariant scalar product on $\Lambda^2\R^{n-k}$ (for
instance the opposite of the Cartan-Killing form
on $\Lambda^2\R^{n-k}=\so(n-k)$). Then a homotopy inverse to $\pi$ is given by
$\sigma\colon U\in \Gr_k(\Lambda^2\R^{n-k}) \mapsto (U,g_0|_U)$.
Clearly, $\pi\circ\sigma =\id_{\Gr_k(\Lambda^2\R^{n-k})}$ and
$H((U,g),t):=(U,(1-t)g+tg_0|_U)$ defines an $\O(n-k)$-equivariant
homotopy $\sigma\circ\pi\simeq\id_{\eS\gamma_{k,\Lambda^2\R^{n-k}}^*} $.
\end{proof}

\begin{remark}\label{aftertheo}
\rm{Using Theorem \ref{general}, one can describe some special cases for
$\NN_{n,k}$.

\noindent $\bullet$ For $k=0$,  $\NN_{n,0}$ is  a point,
corresponding to the $n$-dimensional abelian Lie algebra.

\smallskip

\noindent $\bullet$ For $k=1$, $\Gr_1(\Lambda^2\R^{n-1})$ is 
homeomorphic to the real projective
space $\RP^{ \binom{n-1}{2}-1}$. Moreover,
$A\in\Lambda^2\R^{n-1} \cong \so(n-1)$ is conjugate to a block-diagonal
matrix with $2\times 2$ blocks
$\left(\begin{array}{cc}
0                      &   \lambda_i  \\
-\lambda_i   &   0
\end{array}\right)
$
on the diagonal and such that
$0\leq\lambda_1\leq\ldots\leq\lambda_{\left[\frac{n-1}{2}\right]}$.
Hence
$\Gr_1(\Lambda^2\R^{n-1})/\O(n-1)\approx\Delta^{\left[\frac{n-1}{2}\right]-1}$
is
homeomorphic to a
$(\left[\frac{n-1}{2}\right]-1)$-simplex
and
$$\NN_{n,1}\approx \Delta^{\left[\frac{n-1}{2}\right]-1}\times\R^+ \ .  $$
For odd $n$ and $\lambda_1 = \ldots =\lambda_{\frac{n-1}{2}} = 1$ we 
recover the
$n$-dimensional Heisenberg algebras $\ch_n \in\NN_{n,1}$.
\smallskip

\noindent $\bullet$  If $k=\binom{n-k}{2}$, then
$\Gr_k(\Lambda^2\R^{n-k})$ is homeomorphic to a point and
$\NN_{n,k} \approx \eS(\Lambda^2\R^{n-k})^*/\O(n-k)$ is a quotient of the cone
$\eS(\Lambda^2\R^{n-k})^*$.}

\end{remark}

\section{Metric $2$-step nilpotent Lie algebras of dimension $\leq 6$}

In this section we study in detail the case of Lie algebras of 
dimension up to $6$.
We denote by $\NN_{*,*}^0$ the subspace of Lie algebras with $\Tr(j^*j)=2 $.
The whole space $\NN_{*}$ is a cone over $\NN_{*}^0$ whose peak is 
the abelian Lie algebra.
Clearly, for $m\leq 2$, $\NN_{m}$ is a point, the abelian Lie algebra.
For $m=3,4$ we get $\NN^{0}_3= \{ \ch_3 \}$ and $\NN^{0}_4=\{ \ch_3\oplus\R\}$.
For $m=5$ we have $\NN_5 =  \NN_{5,0} \cup \NN_{5,1} \cup \NN_{5,2}$. By
remark \ref{aftertheo},
$\NN_{5,1}^0$ is homeomorphic to an interval with endpoints the Lie algebras
$\ch_3\oplus\R^2$
and $\ch_5$. Let now $\cn_5\in\NN_{5,2}$ denote a Lie algebra with
isometric $j\colon\R^2\to\so(3)$; all such Lie algebras are 
isometrically isomorphic.
We will se later that the closure $\overline{\NN_{5,2}^0}$ is 
homeomorphic to an
interval with endpoints $\ch_3\oplus\R^2$ and $\cn_5$.
For any $m\leq n$ there are embeddings
$\NN_{m,k}\hookrightarrow\NN_{n,k}$, $\cn\mapsto\cn\oplus\R^{n-m}$. Thus all
the spaces of Lie algebras above appear in $\NN_6$.

In the sequel we will show that $\NN_6$ is a cone over a
contractible $4$-dimensional simplicial complex pictured in Figure~1.
The decomposition \eqref{decomp} becomes
$$\NN_6=\NN_{6,0} \cup \NN_{6,1} \cup \NN_{6,2} \cup \NN_{6,3} \ . $$
 From remark \ref{aftertheo} we have
\begin{eqnarray*}
\NN_{6,0} &\approx & *	\ , \\
\NN_{6,1} &\approx &  [0,1] \times\R^+	\ , \\
\NN_{6,3} &\approx & \eS(\Lambda^2\R^{3})^*/\O(3)      \ .
\end{eqnarray*}

\subsection{Invariants for $\NN_6$}

The subsequent simultaneous description of the pieces of 
\eqref{decomp} and their glueing
relies on the isomorphism of Lie algebras
\begin{equation}\label{SO4split}
\so(4)\cong\so(3)\oplus\so(3)=\su_+(2)\oplus\su_-(2)=\R^3_+\oplus\R^3_- \ .
\end{equation}
Under the identification \eqref{SO4split} the action of $\SO(4)$ on $\so(4)$
translates to the (dual of the) usual action of $\SO(3)\times\SO(3)$
on $\R^3_+\oplus\R^3_-$. The whole orthogonal group in addition contains an
element $\tau\in\O(4)$ of determinant $-1$ which interchanges the factors.
Explicitely, the isomorphism \eqref{SO4split} is given by mapping
$$\xi e_1^\pm+    \psi  e_2^\pm  + \chi e_3^\pm =
\left(\begin{array}{cccc}
i \xi   & \psi +i \chi  \\
-\psi +i \chi & - i \xi
\end{array}\right) \in\su_\pm(2)
$$
to
$$
\left(\begin{array}{cccc}
0      & \xi   & \psi    & \chi  \\
-\xi   & 0     & -\chi  & \psi   \\
-\psi  & \chi & 0       & -\xi    \\
-\chi  & -\psi  & \xi    & 0      \\
\end{array}\right)
  \quad  {\mbox {and}} \quad
\left(\begin{array}{cccc}
0      & \xi   & \psi    & \chi  \\
-\xi   & 0     & \chi  & -\psi   \\
-\psi  & -\chi & 0       & \xi    \\
-\chi  & \psi  & -\xi    & 0      \\
\end{array}\right)
$$
  for "$-$" and "$+$" respectively. The diagonal matrix 
$\mathrm{diag}(-1,1,1,1)$ acts as the
  involution $\tau$.
We denote the two components of $j$ by $j_\pm\colon\R^2\to \R^3_\pm$.
The spectra of $j_\pm^*j_\pm$ and the trace of $j_-^*j_- j_+^*j_+$
are invariant under the $\O(2)\times\O(4)$-action, up to interchanging $\pm$.

We claim that these data  suffice to determine the equivalence
class of $j$ under the $\O(2)\times\O(4)$-action:

Clearly, the entire matrices $j_-^*j_-$ and $j_+^*j_+$ determine
$j$ up to the action of $\O(4)$.
If both $j_-^*j_-$ and $j_+^*j_+$ have two identical eigenvalues, 
then both matrices are diagonal
for any orthonormal basis of $\R^2$.
Otherwise, after possibly using $\tau$ to permute $\pm$, we may assume that
$j_-^*j_-$ has two different eigenvalues $\alpha_-$, $\beta_-$ and that
$e_1,e_2$ are the respective eigenvectors.
If $j_+^*j_+ =
\left(\begin{array}{cc}
x& z \\
z& y
\end{array}\right),$
then $\Tr(j_+^*j_+) = x+y $, $\det(j_+^*j_+)=xy-z^2$ and
$\Tr(j_+^*j_+j_-^*j_-) = \alpha_- x + \beta_- y $
determine $x,y\geq 0$ and $z$ up to sign. Since the sign of $z$ can 
be changed by
conjugation with $\left(\begin{array}{cc}-1& 0 \\0 & 1 
\end{array}\right),$ we may assume $z\geq 0$.
Thus, all of $j_+^*j_+$ is determined by the above invariants.
\medskip

Let $\Spec(j_\pm^*j_\pm)=\{\alpha_\pm, \beta_\pm\}$, with 
$0\leq\alpha_\pm\leq\beta_\pm$ and
$\Tr(j_+^*j_+j_-^*j_-)=t$.
The possible range for $t$ in dependence of $\alpha_\pm, \beta_\pm$ is obtained
by solving
$$t=\alpha_- x + \beta_- y \ ,\quad  x+y = \alpha_+ + \beta_+  \ 
,\quad xy - z^2 = \alpha_+ \beta_+ \ ,\quad z \geq 0 \ . $$
We get
$$t\in I_{\alpha_\pm, \beta_\pm}=
\left[ \alpha_- \beta_+ + \alpha_+\beta_-     ,  \alpha_- \alpha_+ + 
\beta_- \beta_+  \right] \ . $$

Let $S$ denote the set of all $5$-tuples
$(\alpha_\pm, \beta_\pm, t)$ satisfying the above conditions and 
subject to the relation induced from $\tau$, i.e.
$$S:=\left. \left\{(\alpha_-, \beta_-,\alpha_+, \beta_+, t) \in \R^5\left|
	\begin{array}{l}
		 0\leq\alpha_-\leq\beta_- , \\
		 0\leq\alpha_+\leq\beta_+ , \\
		t\in I_{\alpha_\pm, \beta_\pm}
	\end{array}
\right.\right\} \right/ \sim$$
with the identification
$$(\alpha_-,\beta_-,\alpha_+,\beta_+,t) \sim (\alpha_+,\beta_+ 
,\alpha_-,\beta_-,t) \ . $$
We have
\begin{theorem}
The closure of $\NN_{6,2}$ is homeomorphic to $S$ under the map
\begin{eqnarray*}
\Psi\colon\overline{\NN_{6,2}} &\to& S \\
j =(j_-,j_+)&\mapsto& (\Spec(j_-^*j_-),\Spec(j_+^*j_+),\Tr(j_+^*j_+j_-^*j_-))
\end{eqnarray*}
\end{theorem}
\begin{proof} We have already shown that the above map is bijective. 
It is continuous since the spectrum of a
matrix depends continuosly on its entries.
Since
$$\Tr(j^*j)=\Tr(j_-^*j_-)+\Tr(j_+^*j_+)=\alpha_-(j) 
+\beta_-(j)+\alpha_+(j)+\beta_+(j),$$
we get  that, for all $ r>0$, $\Psi$ defines a continuous bijection
$$
\begin{array}{l}
  \{j\mid \Tr(j^*j)\leq r\}/_{\O(2)\times\O(4)} \leftrightarrow\\
S\cap\{(\alpha_-,\beta_-,\alpha_+,\beta_+, t)\in\R^5\mid \alpha_- 
+\beta_-+\alpha_+ +\beta_+\leq r\}/_{\sim}
\end{array}
$$
which is a homeomorphism since these sets are compact. It follows 
that $\Psi$ is a homeomorphism.
\end{proof}

The spaces $\NN_{6,1}$ and $\NN_{6,3}$ are treated similarly.
For $\NN_{6,1}$ we have to deal with maps $j\colon\R\to\so(5)$.
Any such map is conjugate to some $j\colon\R\to\so(4)\subset\so(5)$. 
Extending $j$ by $0$ to
a map $\R^2 \to\so(4)$, we can identify $\NN_{6,1}$ with a subset of 
$\partial\overline{\NN_{6,2}}$.
In the terminology above, both components
$j_+$ and $j_-$  have only one nonvanishing eigenvalue, 
$0\leq\beta_+\leq\beta_-$ respectively. Moreover
$\Tr(j_+^*j_+j_-^*j_-)= \beta_- \beta_+$ gives no new invariant on $\NN_{6,1}$.

For $\NN_{6,3}$, we observe that the imbedding 
$\so(3)\hookrightarrow\so(4)$ induced from
$\R^3\hookrightarrow\R^4$ translates under \eqref{SO4split} to the 
skew-diagonal map
$\so(3)\ni X \mapsto \frac 1 2 (X,-X) \in \so(3)\oplus\so(3)=\so(4)$. 
Thus, for $j=(j_-,j_+)\in\NN_{6,3}$
we have $j_+=-j_-$. Hence,  $\Spec(j_-^*j_-) = \Spec(j_+^*j_+) = 
\{\omega,\alpha,\beta\}$ with
$0\leq\omega\leq\alpha=\alpha_-=\alpha_+\leq\beta=\beta_-=\beta_+$ 
and $t=\omega^2+\alpha^2+\beta^2$.

As a whole,  $\NN_6$ is a cone over the set $\NN_6^0$
of those isometric isomorphism classes with
$\Tr(j^*j)=\Tr(j_-^*j_-) + \Tr(j_+^*j_+)=\alpha_- + \beta_- + 
\alpha_+ +\beta_+ =2$. Its
peak is the abelian Lie algebra $\R^6$. The following picture 
illustrates the set $\NN_6^0$,
where the invariant $t$ is omitted over the interior of $\NN_{6,2}$.
We have chosen a fundamental domain for the $\tau$-action such
that the parameters $\alpha_\pm, \beta_\pm, t$ always satisfy
$\beta_- -\alpha_- \geq \beta_+ -\alpha_+ $. We then only need to identify
$(\alpha_\pm, \beta_\pm, t)\sim (\alpha_\mp, \beta_\mp, t)$ if 
$\beta_- -\alpha_- = \beta_+ -\alpha_+$.
On the right hand face of $\NN_{6,2}$ this requires to identify
the two triangles by the reflection indicated by the two arrows 
$\leftrightarrow$.
The dots $\bullet$ mark standard representatives for the seven 
different isomorphism classes
of Lie algebras and are identified in the next section.

\def\pkt{\circle*{6}}
\def\ccx{\color[rgb]{0,0.7,0}}	
\def\ccy{\color[rgb]{0.2,0.2,0.5}}	
\def\ccaa{\color[rgb]{0,0,1}}	
\def\ccbb{\color[rgb]{1,0,0}}	
\def\ccidf{\color[rgb]{0.6,0.6,0}}	

\begin{figure}[h!]
\begin{picture}(250,260)(60,-10)
\unitlength94000sp
\put(100,0){\pkt\put(-4,-9){\footnotesize $(1,1,0,0,0)\cong\ch_3^\C$}}
\put(160,90){	\put(36,24){\vector(-3,-2){34}}
	\put(38,24){\footnotesize $(0,0,1,1,0)\cong\ch_3^\C$}
	\pkt
	}

\put(100,130){	\put(0,0){\line(3,-2){60}}
	\put(0,0){\line(-3,-2){60}}
	\put(-70,15){\parbox[c]{3cm}{\footnotesize $(0,1,0,1,0)$ 
\\$\cong \ch_3\oplus \ch_3$}}
	\put(-30,15){\vector(2,-1){28}}

	\put(-2,0){\ccx\line(1,3){10}}\put(-1,0){\ccx\line(1,3){10}}\put(0,0){\ccx\line(1,3){10}}

	\put(-6,-4){\ccx\line(1,3){9}}\put(-12,-8){\ccx\line(1,3){8}}\put(-18,-12){\ccx\line(1,3){7}}
	\put(-24,-16){\ccx\line(1,3){6}}\put(-30,-20){\ccx\line(1,3){5}}\put(-36,-24){\ccx\line(1,3){4}}
	\put(-42,-28){\ccx\line(1,3){3}}\put(-48,-32){\ccx\line(1,3){2}}

	\put(3,-9){\ccx\line(1,3){8}}\put(6,-18){\ccx\line(1,3){7}}\put(9,-27){\ccx\line(1,3){6}}

\put(12,-36){\ccx\line(1,3){6}}\put(15,-45){\ccx\line(1,3){5}}
\put(18,-54){\ccx\line(1,3){4}}
\put(21,-62){\ccx\line(1,3){3}}
\put(24,-74){\ccx\line(1,3){3}}

	\put(1,0){\line(1,-3){28}}
	\pkt
	}

\put(207,92){\ccbb\put(0,0){\line(-3,2){100}}
	\put(0,0){\line(-5,-3){78}}
		\put(-15,10){\line(-5,-3){65}}\put(-30,20){\line(-5,-3){53}}
		\put(-45,30){\line(-5,-3){42}}
		\put(-60,40){\line(-5,-3){30}}\put(-75,50){\line(-5,-3){18}}
	\put(0,-30){$\NN_{6,3}$}\put(0,-24){\vector(-2,1){25}}
	\pkt\put(0,-4){\parbox[c]{3cm}{\footnotesize 
$\omega=\alpha=\beta=1/3$, \\ $t=1/3 $}}
	}

\put(100,65){\put(-105,-20){\footnotesize $(1/2,1,0,1/2,1/4)$}
	\put(-45,-15){\vector(3,1){42}}
	\pkt
	}
\put(100,0){\line(2,3){60}}
\put(100,0){\line(-2,3){60}}
\put(100,0){\line(0,1){130}}
\put(130,110){\put(38,24){\footnotesize $(0,1/2,1/2,1,1/4)$}
	\put(36,24){\vector(-3,-2){34}}
	\pkt}

\put(30,5){\put(20,6){\vector(1,1){40}} $\NN_{6,2}$}


\put(120,75){\ccidf\put(-4,-6){\vector(2,3){12}}\put(4,6){\vector(-2,-3){12}} 
}

\put(40,90){\line(1,0){120}}

\put(130,45){\vector(-2,-3){10}}\put(130,45){\vector(-2,-3){20}}
\put(130,45){\vector(2,3){10}}\put(130,45){\vector(2,3){20}}
\put(100,130){\vector(0,-1){32}}\put(100,130){\vector(0,-1){100}}
\put(100,130){\vector(3,-2){15}}\put(100,130){\vector(3,-2){45}}

\put(108,158){\ccy
	\put(0,0){\line(1,-5){23}}
		\put(49,-7){\vector(-3,-2){41}}
		\put(51,-10){$\NN_{5,2}$ }
	}
\put(130,45){\ccy\pkt\put(-4,-9){\footnotesize 
$(1/2,1/2,1/2,1/2,1/2)\cong\cn_5\oplus\R$}}

\put(40,90){\ccaa\pkt
	\put(-50,-4){\parbox[c]{3cm}{\footnotesize $(0,2,0,0,0)$ 
\\$\cong \ch_5\oplus\R$}}
	\put(-54,30){$\NN_{5,1}=\NN_{6,1}$ \put(0,3){\vector(1,0){38}}}
	}
\put(108,158){\ccaa
	\put(0,0){\line(-1,-1){68}}\put(1,-1){\line(-1,-1){68}}\put(-1,1){\line(-1,-1){68}}
	\put(0,0){\line(-1,-1){68}}\put(0,-1){\line(-1,-1){68}}\put(-1,0){\line(-1,-1){68}}
	\put(-4,6){\footnotesize $(0,1,0,1,1)\cong \ch_3\oplus\R^3$}
	\pkt
	}

\end{picture}
\label{bigfig}
\caption{$\NN_6^0$}
\end{figure}

\subsection{Classification of $6$-dimensional $2$-step nilpotent Lie algebras}

Next, we  determine the isomorphism classes (disregarding the metric) 
of $2$-step nilpotent
$6$-dimensional Lie algebras and exhibit canonical representatives. 
To this end, we compute the
action of $\GL(2,\R)\times\GL(4,\R)$ on our invariants.

\begin{remark}
{\rm Using the $\GL(2,\R)$-action for any $j$ we find an equivalent 
one which is an isometric monomorphism.}
\end{remark}
\begin{remark}
{\rm For isometric $j$, we can simultaneously diagonalize $j_+^*j_+$ 
and $j_-^*j_-$. For $\NN_{6,2}$, this yields the relations
$\alpha_- + \beta_+ =1=\alpha_+ + \beta_-$ and $t=\alpha_-\beta_+ + 
\alpha_+  \beta_- $. For $\NN_{6,3}$, we have
$\omega=\alpha=\beta=\frac 1 2$ and for $\NN_{6,1}$ we get 
$\alpha_-=\alpha_+=0$, $\beta_+ + \beta_- =1$  and $t=\beta_+
\beta_-$.}
\end{remark}
A pair $(T,S)\in\GL(2,\R)\times\GL(4,\R)$ acts on $j$ by replacing $j(z)$ with
$S^*j(Tz)S$.
In the bases $e_i^\pm$ for $\R^3_\pm=\so_\pm(3)\subset\so(4)$ and 
$(e_1, e_2)$ for $\R^2$, we
can put $j$ in the form
\begin{equation}\label{jform}
j=\left(\begin{array}{cc}
a_-  & 0    \\
0     & b_- \\
0     & 0      \\
p      & r      \\
0       & q \\
0      & 0
\end{array}\right)
\end{equation}
with $0\leq a_- \leq b_- $, $0\leq p,q,r$, using only the 
$\O(2)\times\O(4)\subset\GL(2,\R)\times\GL(4,\R)$
action. The coefficients $ a_- , b_- , p,q,r \in\R_0^+$ are 
determined from the invariants by solving the equations
$$ \begin{array} {l} a_-^2=\alpha_-, \quad b_-^2=\beta_-, \quad 
p^2+q^2+r^2=\alpha_+ + \beta_+,\\
  \quad p^2q^2=\alpha_+\beta_+,
\quad \alpha_- p^2 + \beta_- (q^2+r^2) = t.
\end{array}$$

We first assume that $j$ is isometric, i.e. $r=0$, 
$a_-^2+p^2=1=b_-^2+q^2$ and $p^2=\beta_+$, $q^2=\alpha_+$.
Possibly interchanging $\pm$ we may also suppose 
$\beta_-\geq\beta_+$. Then, the only free invariants
for this case are $\alpha_\pm$ and satisfy the conditions 
$\alpha_-\geq\alpha_+$ and $\alpha_+\leq 1-\alpha_-$.

By means of the isomorphism \eqref{SO4split} $j$ defines the 
homomorphism $j\colon\R^2\to\so(4)$ given by
$$
j(u,v)  = \left(\begin{array}{cccc}
0	& (a_-+p)u 	& (b_-+q)v& 0 		\\
-(a_-+p)u 	& 0	& 0	& (b_--q)v 	\\
(-b_- -q)v	& 0 	&  0 	& (-a_-+p)u  		\\
0 	&(q-b_-)v	&(a_- -p)u 	& 0 		\\
\end{array}\right)
$$
In case $b_- - q \neq 0$, this is equivalent to the matrix with coefficients
$(a_- , b_-^\prime , p, 0, 0)$, $b_-^\prime=\sqrt{(b_-+q)(b_--q)}$, 
via the matrix
$$
S:=\left(\begin{array}{cccc}
\lambda^{-1} & 0 & 0 & 0 \\
0 & \lambda & 0 & 0 \\
0 & 0 & \lambda^{-1} & 0  \\
0 & 0 & 0 & \lambda \\
\end{array}\right)
$$
with
$\lambda=\displaystyle\left(\frac{b_-+q}{b_--q}\right)^{1/4}$.
By rescaling $v$ we can keep $j$ isometric.
Similarly, in case $a_- - p > 0$, we use a matrix
$$
T:=\left(\begin{array}{cccc}
\lambda^{-1} & 0 & 0 & 0 \\
0 & \lambda^{-1} & 0 & 0 \\
0 & 0 & \lambda & 0  \\
0 & 0 & 0 & \lambda \\
\end{array}\right)
$$
with
$\lambda=\displaystyle \left(\frac{a_-+p}{a_--p}\right)^{1/4}$
to see that any $j$ with coefficients $(a_- , b_- , p, q ,0)$ is 
equivalent to one with
coefficients $(a_-^\prime , b_- , 0, q ,0)$ where \linebreak
$a_-^\prime=\sqrt{(a_-+p)(a_--p)}$. \newline
In case $a_- - p < 0$ we can replace $(a_- , b_- , p, q ,0)$ by $(0 , 
b_- , p, q ,0)$ by means of the above matrix
$T$ with $\lambda=\left(\frac{a_-+p}{p - a_-}\right)^{1/4}$. In order 
to keep $j$ isometric, we rescale $u$.

The diagram below visualizes the subset of $\NN_{6,2}$ represented by 
isometric $j$.
With respect to the action of $\GL(2,\R)\times\GL(4,\R)$ it 
decomposes into four isomorphism classes
of Lie algebras indicated by the components in the picture.

\begin{figure}[h!]
\begin{picture}(200,140)(0,-20)
\unitlength1pt
\put(-15,-10){\scriptsize $\alpha_-=\alpha_+=0$}
\put(85,110){\scriptsize $\alpha_-=\alpha_+=1/2$}
\put(85,-10){\scriptsize $\alpha_-=1/2$, $\alpha_+=0$}
\put(185,-10){\scriptsize $\alpha_-=1$, $\alpha_+=0$}
\put(0,0){\pkt}
\put(100,0){\pkt}
\put(200,0){\pkt}
\put(100,100){\pkt}
\put(0,0){\line(1,0){95}}
\put(95,0){\line(0,1){95}}
\put(0,0){\line(1,1){95}}
\put(100,0){\line(0,1){95}}
\put(200,0){\line 
(-1,0){95}}
\put(105,0){\line (0,1){95}}
\put(200,0){\line 
(-1,1){95}}
\end{picture}
\label{fig}
\caption{$\DD_{6,2}$}
\end{figure}

\begin{remark} 
\label{rem3}
{\rm If $j$ is not isometric and has invariants 
$(\alpha_\mp,\beta_\mp,t)$ we first compute the coefficients
$( a_-, 
b_- , p,q,r) $ to write $j$ in the shape \eqref{jform}. With 

$$
A:=\left(\begin{array}{cc}
1 & 
\displaystyle\frac{-pr}{(a_-^2+p^2)}  \\
0 &1
\end{array}\right) \ 
.
$$
and $B=\mathrm{diag}(1/\|jAe_1\|,1/\|jAe_2\|)$ we get that $jAB$ 
is isometric.
Computing $\alpha_\pm(jAB)$, the isomorphism type of 
$j$ can be determined.}
\end{remark}

In $\NN_{6}$ we get the 
following isomorphism types,
where the parameters are given for 
isometric $j$. 
\begin{enumerate}
\item Any Lie algebra in 
$\NN_{6,1}$ is isomorphic to a Lie algebra with 
parameters
$\alpha_\pm=0=\beta_+ , \ \beta_-=1 , \ t= 0$ or 
$\alpha_\pm=0, \beta_+= \beta_-=1/2 , \ t= 1/4$.
The first type is 
the product $(0,0,0,0,0, 12 + 34)=\R^3\oplus \ch_3$.
The second type 
is a product $(0,0,0,0,0, 12)=\ch_5\oplus \R$. 
\item In $\NN_{6,2}$ 
we have four isomorphism types corresponding to \newline
	(a) 
$\alpha_\pm < 1/2$, which gives $(0,0,0,0, 12, 34)=\ch_3\oplus 
\ch_3$\newline
(b) $\alpha_+ < 1/2=\alpha_-$, $(0,0,0,0, 12, 14 + 
23)$ \newline
	(c) $\alpha_+ = 1/2 = \alpha_-$,
		$(0, 
0,0,0, 12, 13) = \cn_5\oplus \R$ with
		$\cn_5\in\NN_{5,2}$ 
the unique non trivial isomorphism type \newline
(d)  $\alpha_+ < 1/2 
< \alpha_-$, $(0,0,0,0, 13 + 42, 14 + 23)=\ch^\C_3$.
	
\item Any 
Lie algebra in $\NN_{6,3}$ is isomorphic to $(0,0,0,12,13, 23)$
with 
$\omega=\alpha=\beta=1/2$.
\end{enumerate}

\begin{figure}[h!]
\begin{picture}(200,130)(0,-20)
\unitlength1pt
\put(-15,-10){ 
$\ch_3\oplus\ch_3$}
\put(85,108){ $\cn_5\oplus 
\R$}
\put(85,-10){\scriptsize $(0,0,0,0,12,14+23)$}
\put(205,-10){ 
$\ch_3^\C$}
\put(55,40){$\ch_3\oplus\ch_3$}
\put(130,40){$\ch_3^\C$}
\put(0,0){\pkt}
\put(100,0){\pkt}
\put(200,0){\pkt}
\put(100,100){\pkt}
\put(0,0){\line(1,0){95}}
\put(95,0){\line(0,1){95}}
\put(0,0){\line(1,1){95}}
\put(100,0){\line(0,1){95}}
\put(200,0){\line 
(-1,0){95}}
\put(105,0){\line (0,1){95}}
\put(200,0){\line 
(-1,1){95}}
\end{picture}
\label{fig2}
\caption{Isomorphism classes 
in $\DD_{6,2}$}
\end{figure}

\begin{remark}\label{rem1}   {\rm A 
Lie algebra $c \in \LLL_n \subset 
\Hom(\Lambda^2\R^n,\R^n)$  is said 
to
degenerate to another Lie algebra $\tilde c$, if $\tilde c$ is 

represented by a structure which lies in the Zariski closure 
of
the
$\GL(n,\R)$-orbit of a structure which represents $c$. In this 
case 
the entire $\GL(n,\R)$-orbit of $\tilde c$ in $\LLL_n$
lies in 
the closure of the orbit of $c$ \cite{bu, l2}.  Recall that 
$c$ 
degenerates to $\tilde c$ if there exist $g_s \in
\GL(n,\R)$ such 
that
$\lim_{s\to 0} g_s \cdot c=\tilde c$. Using this, it is easy to 
see 
that the  Lie algebras $\ch_3^{\C}, \ch_3 \oplus \ch_3$, 
$(0,0,0,0, 12, 14 + 23)$
 all degenerate to $\N_5\oplus\R$ (the top 
point  in Figure~3).}
\end{remark}

\begin{remark}  {\rm Using Remark 
\ref{rem3}, one can determine the structure equations for any 
6-dimensional 2-step nilpotent
Lie algebra. As an example, from the 
isomorphism
\begin{equation} \label{isoforms}
{\mathfrak {so}} (4) 
\cong \Lambda^2 \R^4
\end{equation}
we will give the structure 
equations for the Lie algebras in  $\DD_{6,2}$.

 Indeed, if one 
fixes a non-zero element $w \in
\Lambda^4 \R^4$, one can consider the 
bilinear form $\phi$ of
signature (3,3) on $\Lambda^2 \R^4$ defined 
by
$\sigma \wedge \tau = \phi (\sigma, \tau) w$.

Given an 
orientation and a metric $g$ on $\R^4$ (and so on $\Lambda^2
\R^4$), 
there is an $SO(4)$-decomposition
\begin{equation} 
\label{decompforms}
\Lambda^2 \R^4 = \Lambda^2_+ \oplus 
\Lambda^2_-,
\end{equation}
where $\Lambda^2_{\pm}$ are the 
eigenspaces of the conformally
invariant involution $*$ of $\Lambda^2 
\R^4$ for which $\phi(*
\sigma, \tau) = g (\sigma, \tau)$. From a 
representation-theoretic
point of view, \eqref{decompforms} is 
equivalent to the Lie
algebra splitting \eqref{SO4split}.

\noindent 
If one chooses a basis $\{ e^1, e^2, e^3, e^4 \}$ of $\R^4$ such 
that
$w = e^1 \wedge e^2 \wedge e^3 \wedge e^4$, 
then
$$
\begin{array} {l}
\Lambda^2_+ = {\rm span \, }\{e^1 \wedge 
e^2 + e^3 \wedge e^4, e^1
\wedge e^3 + e^4 \wedge e^2, e^1 \wedge e^4 
+ e^2 \wedge
e^3\},\\
\Lambda^2_- = {\rm span \, }\{e^1 \wedge e^2 - 
e^3 \wedge e^4, e^1
\wedge e^3 - e^4 \wedge e^2, e^1 \wedge e^4 - e^2 
\wedge
e^3\} .
\end{array}
$$
Using \eqref{isoforms} and the 
embeddings \eqref{SO4split} one has
the following 
identifications
$$
\begin{array} {l}
e_1^+ \sim e^1 \wedge e^2 + e^3 
\wedge e^4,\\
e_2^+ \sim e^1 \wedge e^3 - e^2 \wedge e^4,\\
e_1^- 
\sim e^1 \wedge e^2 - e^3 \wedge e^4,\\
e_2^- \sim e^1 \wedge e^3 + 
e^2 \wedge e^4
\end{array}
$$
and thus $\N_{(\alpha_+,\alpha_-)}$ has 
structure equations
\begin{equation}\label{str-eq}
\begin{array} 
{l}
d e^i = 0,\qquad i=1, \dots , 4 \ ,\\
d e^5 = (a_- + p) \, e^1 
\wedge e^2 + (a_- -
p) \, e^3 \wedge e^4,\\
d e^6 = (b_- + q) \, e^1 
\wedge e^3 - (b_- -
q) \, e^2 \wedge 
e^4,
\end{array}
\end{equation}
where 
$$
a_-= \sqrt{\alpha_-}, \quad 
b_- = \sqrt{\beta_-}, \quad  p= \sqrt{\beta_+}, \quad q = 
\sqrt{\alpha_+}.
$$
}
\end{remark}

\begin{remark}

\rm{ Next, we 
compute the infinitesimal rank of a Lie algebra in $\DD_{6,2}$.

The 
rank of a geodesic in a Riemannian manifold $M$ is the dimension of 
the
real vector space of parallel Jacobi fields along it. The rank 

$\mathrm{rk}(M)$ of $M$
is the minimum of the ranks of all its 
geodesics.
Recall that the Jacobi-operator $R_v$ is the endomorphism 
of $T_p M$
given by $w \mapsto R_{v,w} v$. The infinitesimal rank 

$\mathrm{inf}\mathrm{rk}(M)$ of $M$ is the
minimal dimension of the 
kernels of its Jacobi-operators \cite{sa2}.
A Riemannian manifold $M$ 
has higher (infinitesimal)
rank if $(\mathrm{inf})\mathrm{rk}(M)\geq 
2$.

First, we use the structure equations \eqref{str-eq}  to 
compute the 
curvature tensor
with respect to the metric $g$ for 
which the forms $(e^i)$ are dual to 
an orthonormal basis 
$(e_i)$.

The  non vanishing components $R_{ijhk}=g(R_{e_i, e_j} e_h, 
e_k)$ of 
the Riemannian curvature tensor are:
$$
\begin{array} 
{l}
R_{1212} = -\frac{3}{4} (a_- + {p})^2,\\
R_{1234} = -\frac{1}{2} 
a_-^2 + \frac{1}{2} {p}^2 + \frac{1}{4} b_-^2 - \frac{1}{4} 
{q}^2 = 
R_{3412},\\
R_{1313} = -\frac{3}{4} (b_- - {q})^2,\\
R_{1324} = 
-\frac{1}{4} a_-^2 +  \frac{1}{4} {p}^2 + \frac{1}{2} b_-^2 - 
\frac{1}{2} 
{q}^2= R_{2413},\\
R_{1423} = \frac{1}{4} a_-^2 - 
\frac{1}{4} {p}^2 + \frac{1}{4} b_-^2 - \frac{1}{4} 
{q}^2 = 
R_{2314},\\
R_{1456} =  \frac{1}{2} {p} {q} - \frac{1}{2} a_- b_- = 
R_{5614},\\
R_{1515} = \frac{1}{4} (a_- + {p})^2 = R_{2525}, 
\\
R_{1546} = -\frac{1}{4} (b_- + {q}) (a_- - {p}) = 
R_{4615},\\
R_{1616} = \frac{1}{4} (b_- + {q})^2 = 
R_{3636},\\
R_{1645} = \frac{1}{4} (b_- - {q}) (a_- + {p}) = 
R_{4516},\\
R_{2356} = -\frac{1}{2} a_- b_- - \frac{1}{2} {p} {q}= 
R_{5623},\\
R_{2424} = -\frac{3}{4} (b_- - {q})^2,\\
R_{2536} = 
-\frac{1}{4} (a_- - {p}) (b_- - {q}) = R_{3625},\\
R_{2626} = 
\frac{1}{4} (b_- - {q})^2 = R_{4646},\\
R_{2635} =  \frac{1}{4} (b_- 
+ {q}) (a_- + {p})=R_{3526},\\
R_{3434} = - \frac{3}{4} (a_- - 
{p})^2,\\
R_{3535} = \frac{1}{4} (a_- - {p})^2 = 
R_{4545}.
\end{array}
$$
The infinitesimal rank is 1 for 
any
$a_-,
{p},
b_-,
q$ (with respect to the above metric), except for 
$(a_-, 
b_-) = (1,0)$ and $(a_-,
b_-) = (\sqrt{2}/2, {\sqrt 2}/2)$. 
Indeed the  Jacobi operator:
$$
R_{e_1 + e_6}: X \mapsto R_{e_1 + 
e_6, X} e_1 + e_6
$$
whose associated matrix is
{\small
$$
\left( 
\begin{array} {cccccc} \frac{1}{4} \eta^2&0&0&0&0&- \frac{1}{4} 
\eta^2\\
0&-\frac{1}{4} [3 \nu^2 - (b_- - q)^2]&0&0&0&0\\
0&0&- 
\frac{1}{2} \eta^2&0&0&0\\
0&0&0& \frac{1}{4} (b_- - 
q)^2&-\rho&0\\
0&0&0&\rho&\frac{1}{4}\nu^2&0\\
- \frac{1}{4} 
\eta^2&0&0&0&0&\frac{1}{4} \eta^2 \end{array} \right)
$$}
(with $\nu 
= a_- + p$, $\eta = b_- + 
q$, $\rho=\frac{1}{4}
[3 p q - 3 a_- b_- 
+
p b_- -
a_- q]$)
  has one dimensional kernel, except for the 
following cases:
\begin{enumerate}
\item $a_- = q$, $b_- = p$;
\item 
$a_- = p = \sqrt{2}/2$.
\end{enumerate}

If one considers, in 
addition,  the  Jacobi operator:
$$
R_{e_2 + e_5}: X \mapsto R_{e_2 + 
e_5, X} e_2 + e_5 \,
$$
its associated matrix 
is
{\footnotesize
$$
\left( \begin{array} {cccccc} -\frac{1} {2} 
\nu^2&0&0&0&0&0\\
0 & \frac{1}{4} \nu^2&0&0&0&-\frac{1}{4} 
\nu^2\\
0&0&\frac{1}{4} (a_- - p)^2&0&-\zeta&0\\
0&0&0&-\frac{1}{4}[3 
(b_- - q)^2 - (a_- - 
p)]&0&0\\
0&-\frac{1}{4} \nu^2&0&0&0& 
\frac{1}{4} \nu^2\\
0&0&-\zeta&0&\frac{1}{4} (b_- - 
q)&0
\end{array}
\right) \ ,
$$}
(with $\zeta=\frac{1}{4}[3 a_- b_- + 
3 p 
q + a_-
q + p b_-]$).
Again, the dimension of the kernel of 
$R_{e_2 + e_5}$ is generically 
one. Both $R_{e_1 + e_6}$ and
$R_{e_2 
+ e_5}$ have kernel of dimension bigger than one if
$b_- =0 =
p$ and 
$a_- = b_- =
\sqrt{2}/2$.
These two cases correspond to the Lie 
algebras $\cn_{(1,1)}\cong 
{\mathfrak h}_3 \oplus {\mathfrak
h}_3$ 
and $\cn_({1/2,1/2)} \cong \N_5\oplus\R$, which 
are both Riemannian 
products of
(infinitesimal) rank one Lie algebras. Thus, their 
(infinitesimal) rank is two.}
\end{remark}

\noindent {\it 
Acknowledgments} We are grateful to  Simon Salamon and Sergio 
Garbiero  for helpful discussions and 
useful
references.


\begin{thebibliography}{BB}

\bibitem{bu} D. 
Burde, Degenerations of nilpotent Lie algebras. 
J. Lie Theory  {\bf 
9}  (1999),  no. 1, 193--202.

\bibitem{df} I. Dotti, A. Fino, 
HyperK\"ahler torsion structures 
invariant by nilpotent Lie groups, 
Classical
Quantum Gravity  {\bf 19} (2002),  no. 3, 
551--562.

\bibitem{df2} I.  Dotti, A.  Fino: Abelian hypercomplex 
8-dimensional
 nilmanifolds,
{\it Ann. Glob. Anal. and  Geom.} {\bf 
18} (2000), 47-59.

\bibitem{e} P. Eberlein, Geometry of $2$-step 
nilpotent groups with
a left invariant metric, Ann. Scient. Ec. Norm. 
Sup., 4, {\bf 27}
(1994), 611-660.

\bibitem{e2} P. Eberlein, The 
moduli space of 2-step nilpotent Lie groups of type (p,q), preprint 
(2002).

\bibitem{e3} P. Eberlein, Geometry of 2-step nilpotent Lie 
groups ", preprint (2003).



\bibitem{gms} C. Gordon, Y. Mao, D. 
Schueth,  Symplectic 
rigidity of geodesic flows on 
two-step
nilmanifolds.  Ann. Sci. \'Ecole Norm. Sup. (4)  {\bf 30} 
(1997), 
no. 4, 417--427.

\bibitem{GK} M. Goze, Y. Khakimdjanov, 
Nilpotent Lie algebras,
Mathematics and its Applications, 361, 
Kluwer
Academic Publishers Group, Dordrecht, 1996.

\bibitem{k} A. 
Kaplan, Riemannian nilmanifolds attached to Clifford
modules, 
Geometriae Dedicata {\bf 11} (1981), 127-136.

\bibitem{k2} A. 
Kaplan, On the geometry of groups of Heisenberg type,
Bull.\ London 
Math.\ Soc.\ {\bf 15} (1983), 35-42.

\bibitem{l} J. Lauret, 
Homogeneous nilmanifolds attached to representations
of compact Lie 
groups, Manuscripta math. {\bf 99} (1999), 287-309.

\bibitem{l2}  J. 
Lauret, Degenerations of Lie algebras and geometry of Lie
groups, 
Differential Geometry and its Applications {\bf 18} (2003), 
177-194

\bibitem{M} L. Magnin,
Sur les alg\`ebres de Lie nilpotentes 
de dimension $\leq 7$,  J. Geom.
Phys. {\bf 3} (1986), no. 1, 
119--144.

\bibitem{Sa} S. Salamon, Complex structures on nilpotent 
Lie algebras,
J. Pure Appl. Algebra {\bf 157} (2001), no. 2-3, 
311-333.

\bibitem{sa2} E. Samiou, 2-step nilpotent Lie groups of 
higher rank,
Manuscripta math. {\bf 107} (2002), 
101-110.


\end{thebibliography}
\end{document}